\documentclass[10pt, reqno]{amsart}
\usepackage[dvipsnames]{xcolor}
\usepackage{graphicx}
\usepackage{amsmath}
\usepackage{amsfonts}
\usepackage{amssymb}
\usepackage{bbm}
\usepackage{amsthm}
\usepackage{tikz}
\usepackage[sc]{mathpazo}
\usepackage[T1]{fontenc}
\usepackage{bm}

\usepackage{enumitem}
    \setenumerate{itemsep=5pt} 
    \setitemize{itemsep=5pt} 

\theoremstyle{plain}
\newtheorem{theorem}{Theorem}

\theoremstyle{definition}

\newcommand{\R}{\mathbb{R}}

\allowdisplaybreaks
\usepackage{hyperref}

\begin{document}

\title{The Reciprocal Schur Inequality}
\author[A.~B\"ottcher]{Albrecht B\"ottcher}
\address{Fakult\"at f\"ur Mathematik, TU Chemnitz, 09107 Chemnitz, Germany}
\email{aboettch@mathematik.tu-chemnitz.de}
\urladdr{\url{https://www-user.tu-chemnitz.de/~aboettch/}}

\author[S.R.~Garcia]{Stephan Ramon Garcia}
\address{Department of Mathematics and Statistics, Pomona College, 610 N. College Ave., Claremont, CA 91711, USA}
\email{stephan.garcia@pomona.edu}
\urladdr{\url{http://pages.pomona.edu/~sg064747}}

\author[M. Mitkovski]{Mishko Mitkovski}
\address{School of Mathematical and Statistical Sciences, Clemson University, Clemson, SC 29634, USA}
\email{mmitkov@clemson.edu}
\urladdr{\url{https://mmitkov.people.clemson.edu/}}

\thanks{SRG partially supported by NSF grant DMS-2054002.}

\begin{abstract}
Schur's inequality states that the sum of three special terms is always nonnegative.
This note is a short review of inequalities for the sum of the reciprocals of these terms and
of extensions of the latter inequalities to an arbitrary number of terms and thus to
higher-order divided differences.
\end{abstract}

\maketitle

\noindent
Schur's inequality says that
\begin{equation}\label{Schur}
x^{s}(x-y)(x-z)+y^{s}(y-z)(y-x)+z^{s}(z-x)(z-y)\geq 0
\end{equation}
for $x,y,z >0$ and arbitrary $s \in \R$, with equality if and only if $x=y=z$. The reciprocal
version of this inequality states that, for distinct $x,y,z >0$,
\begin{equation}\label{RecSchur}
\frac{1}{x^{s}(x-y)(x-z)}+\frac{1}{y^{s}(y-z)(y-x)}+\frac{1}{z^{s}(z-x)(z-y)}
\end{equation}
is positive if $s >0$ or $s <-1$, negative if $-1 < s <0$, and zero if $s=0$ or $s=-1$.
It is well known that extension of Schur's inequality to  arbitrarily many variables
$x_1, \ldots, x_n$ is impossible without imposing additional constraints on the variables.
See, e.g.,~\cite{Mit, Wu}.
However, the reciprocal Schur inequality may be extended to an arbitrary number of
variables without any additional requirements. Here it is.

\begin{theorem}\label{Theo1}
Let $x_1, \ldots, x_n$ ($n \ge 2$) be distinct positive real numbers and let $s \in \R$. Then
the sum
\begin{equation}\label{sum}
\sum_{j=1}^n \frac{1}{x_j^s \prod_{k\neq j}(x_j-x_k)}
\end{equation}
is zero if and only if $s \in \{0,-1,-2, \ldots, -(n-2)\}$
and otherwise the sign of this sum
equals the sign of $(-1)^{n+1}s(s+1)\cdots (s+n-2)$.
\end{theorem}

Thus, if, for example, $n=4$, then the sum~(\ref{sum}) is zero for $s\in \{0,-1,-2\}$,
it is positive if $s \in (-\infty,-2) \cup (-1,0)$, and it is negative if $s \in (-2,-1) \cup (1,\infty)$.

\medskip
Theorem~\ref{Theo1} is actually a special case of the following theorem.

\begin{theorem}\label{Theo2}
Let $-\infty \le \alpha < \beta \le \infty$, let $n \ge 2$, and let $f:(\alpha, \beta) \to \R$ be $n-1$
times continuously differentiable in $(\alpha,\beta)$.
Then
\begin{equation}\label{sum3}
\sum_{j=1}^n \frac{f(x_j)}{\prod_{k \neq j}(x_j-x_k)}
\end{equation}
is nonnegative (resp. nonpositive) for arbitrary distinct $x_1, \ldots, x_n$ in $(\alpha, \beta)$
if and only if $f^{(n-1)}$ is nonnegative (resp. nonpositive)
on $(\alpha, \beta)$.
\end{theorem}

Clearly, Theorem~\ref{Theo1} follows from taking $(\alpha,\beta)=(0,\infty)$ and $f(t)=1/t^s$, in which case
\[f^{(n-1)}(t)=(-1)^{n+1}s(s+1)\cdots (s+n-2)\,\frac{1}{t^{s+n-1}}.\]

The sum~(\ref{sum3}) is known as the $n$th divided difference of the function $f$ and is usually denoted by
$f[x_1, \ldots,x_n]$. In~\cite{BGON} and~\cite{FaZwi}, Theorem~\ref{Theo2} is proved as follows.
The Curry--Schoenberg B-spline associated with $x_1 <  \cdots < x_n$ is the function
\begin{equation*}
F(t;x_1,\ldots,x_n) = \frac{n-1}{2} \sum_{j=1}^n \frac{ |x_j - t| (x_j - t)^{n-3}}{ \prod_{k \neq j} (x_j - x_k)}.
\end{equation*}
This is actually a probability density supported on $[x_1,x_n]$. We have in particular, $F(t;x_1,\ldots,x_n) \ge 0$
for $t \in [x_1,x_n]$.
An identity known as Peano's formula says that if $f:[x_1,x_n] \to \R$ is in $C^{n-1}$, then
\begin{equation} \label{P1}
f[x_1,\ldots,x_n]=\frac{1}{(n-1)!} \int_{x_1}^{x_n} f^{(n-1)}(t) F(t;x_1,\ldots,x_n)\,dt.
\end{equation}
See, for instance, \cite[Theorem 3.7.1]{Dav}.
This formula implies the ``if'' portion of the theorem. As for
the ``only if'' part, notice that if $f^{(n-1)}(t_0) < 0$ for some $t_0$, then, by the already proved ``if'' part,
$f[x_1, \ldots,x_n] <0$ for all
$x_1, \ldots,x_n$ sufficiently close to $t_0$.

\medskip
Requiring that $f$ be $n-1$ times continuously differentiable is actually too much. The ultimate
answer to the question about minimal conditions needed for a result like Theorem~\ref{Theo2} was given by Eberhard Hopf
in his dissertation~\cite{Hopf}, defended in 1926. Incidentally, the referees of the dissertation were Erhard Schmidt
and Issai Schur.

Hopf first notes that Theorem~\ref{Theo2} is true if only the existence of $f^{(n-1)}$ on $(\alpha,\beta)$
is required. He refers to a theorem by H.~A.~Schwarz~\cite{Schwarz}, according to which
for an $n-1$ times differentiable function $f$ there is a $t \in (x_1,x_n)$ such that
\begin{equation} \label{HS}
f[x_1, \dots,x_n]=\frac{f^{(n-1)}(t)}{(n-1)!}.
\end{equation}
See also equality (1.33) in \cite{Phil}. Obviously, this formula shows that $f[x_1, \ldots, x_n] \ge 0$
if $f^{(n-1)}(t)\ge 0$ on $(\alpha, \beta)$. To tackle the case where $f^{(n-1)}(t_0) < 0$ for some $t_0$
and continuity of $f^{(n-1)}$ at $t_0$ is not guaranteed,
Hopf has recourse to a result by T.~J.~Stieltjes~\cite{Stie}, which says that $f[x_1, \ldots, x_n]$
converges to $f^{(n-1)}(t_0)/n!$ whenever $x_1, \ldots, x_n$ converge to $t_0$ with the additional
condition that $t_0$ stays between the minimum and the maximum of $x_1, \ldots,x_n$.

Here is Hopf's final theorem on positivity of divided differences.

\begin{theorem}\label{Theo3}
Let $-\infty \le \alpha < \beta \le \infty$, let $n \ge 3$, and let $f:(\alpha, \beta) \to \R$ be
a function.
Then
\begin{equation*}
\sum_{j=1}^n \frac{f(x_j)}{\prod_{k \neq j}(x_j-x_k)}
\end{equation*}
is nonnegative (resp. nonpositive) for arbitrary distinct $x_1, \ldots, x_n$ in $(\alpha, \beta)$
if and only if $f$ is $n-3$ times differentiable and $f^{(n-3)}$ is convex (resp. concave)
on $(\alpha, \beta)$.
\end{theorem}

This is Satz 1 on page 24 of \cite{Hopf}. Repeating Hopf's full proof here is beyond the scope of this note.
We therefore confine ourselves to its basic steps.

For $n=3$ and $x_1 <x_2 <x_3$, the inequality $f[x_1,x_2,x_3]\ge 0$ reads
\[f(x_2) \le f(x_1)\frac{x_3-x_2}{x_3-x_1}+f(x_3)\frac{x_2-x_1}{x_3-x_1},\]
and this holds for all $x_1 <x_2 <x_3$ if and only if $f$ is convex.
So let $n \ge 4$ and suppose the theorem is true for $n-1$.
Without loss of generality assume that $\alpha$ and $\beta$ are finite and that $x_1 < \cdots < x_n$.
Let $f[x_1, \ldots,x_n] \ge 0$ for all $x_1, \ldots, x_n$. The first goal is to prove that then $f$ is differentiable.
Take $0<\delta < (\beta-\alpha)/2$ and fix $n-1$ points $a_1, \ldots, a_{n-1}$ in $(\alpha, \alpha+\delta)$
as well as $n-1$ points $b_1, \ldots, b_{n-1}$ in $(\beta-\delta,\beta)$. A simple identity for divided
differences gives
\begin{align*}
& \frac{n-1}{n}\Big\{f[x_1, \ldots,x_{n-1}]-f[a_1, \ldots,a_{n-1}]\Big\}\\
& = (x_1-a_1)f[a_1,x_1, \ldots,x_{n-1}]+\cdots+(x_{n-1}-a_{n-1})f[a_1, \ldots,a_{n-1},x_{n-1}]
\end{align*}
for every choice of $x_1, \ldots, x_{n-1}$ in $(\alpha+\delta, \beta-\delta)$. As the right-hand side is nonnegative,
it follows that $f[x_1, \ldots,x_{n-1}] \ge f[a_1, \ldots, a_{n-1}]$. It can be shown in an analogous fashion that
$f[x_1, \ldots,x_{n-1}] \le f[b_1, \ldots, b_{n-1}]$. Consequently, there is a constant $M$ such that $|f[x_1, \ldots,x_{n-1}]| \le M$
for all $x_1, \ldots, x_n$. A Hilfssatz proved on page 12 says that this implies that $f$ is differentiable on $(\alpha+\delta, \beta-\delta)$,
and as $\delta >0$ can be made arbitrarily small, one gets differentiability on all of $(\alpha,\beta)$.

Now the induction step. It is based on two remarkable theorems.
The first of them is a generalization of formula~(\ref{HS}). This generalization, Satz 1 on page 9, states
that if $f$ is differentiable, then there are $t_1, \ldots, t_{n-1}$ with $t_j \in (x_j,x_{j+1})$ for all $j$
such that
\begin{equation}\label{Hopf1}
f[x_1, \ldots,x_n]=\frac{1}{n}f'[t_1, \ldots, t_{n-1}].
\end{equation}
The second of the two theorems, Satz 2 on page 11,
says that if $f$ is differentiable, then there is a $t \in (x_1,x_n)$ such that
\begin{equation}\label{Hopf2}
\frac{1}{n}f'[x_1, \ldots,x_{n-1}]=f[t, x_1, \ldots,x_{n-1}].
\end{equation}
Using equalities (\ref{Hopf1}) and (\ref{Hopf2}) and the differentiability of $f$ proved in the previous paragraph,
one gets that $f[x_1. \ldots,x_n] \ge 0$ for all
$x_1, \ldots, x_n$ if and only if $f$ is differentiable and $f'[x_1, \ldots,x_{n-1}] \ge 0$
for all $x_1, \ldots, x_{n-1}$. By the induction hypothesis, the latter is equivalent to
the requirement that the $(n-4)$th derivative of $f'$ exists and is convex, which is the same
as requiring that $f^{(n-3)}$ exists and is convex. This completes Hopf's proof.

\medskip
One of the marvels
in divided differences is that $f[x_1, \ldots,x_n]$ is a complete homogeneous symmetric (CHS) polynomial
if $f$ is a monomial of sufficiently large degree; see, e.g.,~\cite[Lemma~4]{Pales} or~\cite[Theorem~1.2.1]{Phil}. To be precise,
if $q$ is a nonnegative integer and $f(t)=t^{q+n-1}$, then
\[f[x_1, \ldots, x_n]=h_q(x_1, \ldots,x_n):=\sum_{1\le j_1 \le j_2 \le \ldots \le j_q \le n} x_{j_1}x_{j_2} \cdots x_{j_q},\]
with the convention that $h_0(x_1, \ldots,x_n):=1$.
This is in fact Jacobi's
bialternant formula, which says that, for
every nonnegative integer $q$,
\begin{equation*}
h_{q}(x_1,x_2,\ldots,x_n) V(x_1,x_2,\ldots,x_n) =
\det
\begin{pmatrix}
1 & x_1 & x_1^2 & \cdots & x_1^{n-2} & x_1^{q+n-1} \\[2pt]
1 & x_2 & x_2^2 & \cdots & x_2^{n-2} & x_2^{q+n-1} \\
\vdots & \vdots & \vdots & \ddots & \vdots & \vdots \\
1 & x_n & x_n^2 & \cdots & x_n^{n-2} & x_n^{q+n-1} \\
\end{pmatrix},
\end{equation*}
where $V(x_1,x_2,\ldots,x_n) = \prod_{1 \leq i < j \leq n} (x_j - x_i)$
is the Vandermonde determinant.

In the case of even $q=2p$, we have
\[(t^{2p+n-1})^{(n-1)}=(2p+n-1)(2p+n-2)\cdots(2p+1)t^{2p} \ge 0\]
on $(\alpha, \beta)=(-\infty,\infty)$, and hence Theorem~\ref{Theo2} implies that
\begin{equation}\label{sumjj}
\sum_{j=1}^n \frac{x_j^{2p+n-1}}{\prod_{k\neq j}(x_j-x_k)} =h_{2p}(x_1, \ldots,x_n)
\end{equation}
is nonnegative for arbitrary
distinct real numbers $x_1, \ldots,x_n$. 
In the discussion of~\cite{Tao} it is shown that
\begin{equation} \label{ta}
h_q(x_1,x_2,\ldots,x_n)=\frac{1}{q!}\mathbb{E}\big((x_1Z_1+x_2Z_2+\cdots+x_nZ_n)^q\big)
\end{equation}
for arbitrary real numbers $x_1, x_2, \ldots, x_n$,
where $Z_1, Z_2, \ldots,Z_n$ are i.i.d. exponentially distributed random variables with parameter $1$.
Clearly, this is another way to see that $h_q(x_1,x_2,\ldots,x_n)\ge 0$ if $q=2p$ is even.
We want to note that~(\ref{ta}) is the case $g(t)=t^q/q!$ of the more general identity 
\[\mathbb{E}\big(g(x_1Z_1+x_2Z_2+\cdots+x_nZ_n)\big)=H[x_1,x_2, \ldots,x_n]\]
where $H(s)=s^{n-2}G(1/s)$ and $G(s)$ is the Laplace transform of $g(t)$. More about
the latter identity will be said elsewhere.

Actually a sharp lower bound for CHS polynomials of even degree is known.
Namely, a famous theorem by Hunter~\cite{Hun} says that
\begin{equation} \label{Hu}
h_{2p}(x_1, \ldots, x_n)\ge \frac{1}{2^p p!}(x_1^2+\cdots+x_n^2)^p
\end{equation}
for arbitrary $(x_1, \ldots, x_n) \in \R^n$. Combining~(\ref{sumjj}) and~(\ref{Hu})
we see that Hunter's inequality is actually the following reciprocal Schur inequality.

\begin{theorem}\label{Theo4}
Let $x_1, \ldots, x_n$ ($n \ge 2$) be distinct real numbers and let $p$ be a nonnegative integer. Then
\begin{equation*} 
\sum_{j=1}^n \frac{x_j^{2p+n-1}}{\prod_{k\neq j}(x_j-x_k)} \ge \frac{1}{2^p p!}(x_1^2+\cdots+x_n^2)^p.
\end{equation*}
\end{theorem}

\medskip
The following theorem of Farwig and Zwick~\cite{FaZwi} provides us with another lower bound for $f[x_1, \ldots,x_n]$.
It is applicable if $f^{(n-1)}$ exists and is convex (and thus, in particular, continuous).

\begin{theorem}\label{Theo5}
Let $x_1, \ldots, x_n$ ($n \ge 2$) be distinct real numbers and let $f:(\alpha,\beta) \to \R$
be a function whose $(n-1)$th derivative exists and is convex. Then, for distinct $x_1, \ldots, x_n$
in $(\alpha,\beta)$,
\begin{equation*} 
\sum_{j=1}^n \frac{f(x_j)}{\prod_{k\neq j}(x_j-x_k)} \ge \frac{1}{(n-1)!}\,f^{(n-1)}\!\left(\frac{x_1+\cdots+x_n}{n}\right).
\end{equation*}
\end{theorem}

For example, if $f(t)=t^6$, then Theorem~\ref{Theo4} (with $p=2$ and $n=3$) and Theorem~\ref{Theo5}
deliver the lower bounds
\begin{equation*} 
\frac{(x_1^2+x_2^2+x_3^2)^2}{8}\;\:\mbox{and}\;\:
\frac{5}{27}(x_1+x_2+x_3)^4,
\end{equation*}
respectively. If $x_1+x_2+x_3=0$, then the first of the bounds is better than the second, but if
$x_1,x_2,x_3$ are very close to one another, then the second of the two bounds is larger than
the first.

The proof of Theorem \ref{Theo5} given in~\cite{FaZwi} starts with the fact that
\[\int_{x_1}^{x_n} t\, F(t;x_1, \ldots,x_n)\,dt=\frac{x_1+\cdots+x_n}{n}=:c\]
for $x_1 < \cdots < x_n$. Since $f^{(n-1)}$ is convex, we may invoke Jensen's inequality
to obtain
\begin{equation*}
f^{(n-1)}(c)=f^{(n-1)}\!\left(\int_{x_1}^{x_n} t\, F(t;x_1, \ldots,x_n)\,dt\right)
\le \int_{x_1}^{x_n} f^{(n-1)}(t)\, F(t;x_1, \ldots,x_n)\,dt,
\end{equation*}
and Peano's formula tells us that the right-hand side of this inequality is nothing but
$(n-1)!\,f[x_1, \ldots,x_n]$. This completes the proof.

Interestingly, Pe\v{c}ari\'{c} and Zwick~\cite{PeZwi} proved that existence and convexity
of $f^{(n-1)}$ implies that $f[x_1, \ldots,x_n]$ is Schur convex. We first want to notice
that if $f^{(n-1)}$ is continuous on $(\alpha,\beta)$, then $f[x_1, \ldots,x_n]$ may be
continuously extended to a symmetric function on all of $(\alpha,\beta)^n$. Schur convexity means that
$f[x_1, \ldots, x_n] \le f[y_1, \ldots, y_n]$ whenever $x \prec y$, where $x \prec y$ in turn means that
$x_1 \ge \cdots \ge x_n$, $y_1 \ge \cdots \ge y_n$,
$\sum_{j=1}^k x_j \le \sum_{j=1}^k y_j$ for $1 \le k \le n-1$, and $\sum_{j=1}^n x_j = \sum_{j=1}^n y_j$.
For example, if $x_1 \ge \cdots \ge x_n$, then
\begin{equation} \label{SchCo}
\left(\frac{x_1+\cdots+x_n}{n}, \ldots, \frac{x_1+\cdots+x_n}{n}\right) \prec (x_1, \ldots,x_n).
\end{equation}
Taking into account that $f[c, \ldots, c]=f^{(n-1)}(c)/(n-1)!$, we conclude from~(\ref{SchCo}) and the
Schur convexity of $f$ that $f^{(n-1)}(c)/(n-1)! \le f[x_1, \ldots,x_n]$ with $c=(x_1+\cdots+x_n)/n$
for arbitrary points $x_1, \ldots, x_n$.
This is the proof of Theorem~\ref{Theo5} presented in~\cite{PeZwi}.

\medskip
Our interest in the topic of this note was aroused by the question on whether
the positivity of one of the exotic expressions in our recent paper~\cite{BGON} is related to
Schur’s inequality~(\ref{Schur}) of math contest fame. We couldn't prove extensions or generalizations of
Schur's original inequality and rather became aware that we entered the business of
inequalities for divided differences. The interpretation of
these inequalities as reciprocal Schur inequalities seemed us worth
exposing them in a note.

We remark that Schur's original inequality~(\ref{Schur}) is known to hold with $x^s,y^s,z^s$ replaced by $f(x),f(y), f(z)$
with a function $f: (\alpha,\beta) \to \R$
and for $\alpha < x \le y \le z < \beta$ if and only if $f$ belongs to the so-called class $Q$ on $(\alpha,\beta)$.
The class $Q$ was introduced in~\cite{GoLe} and is defined as the set of all functions $f: (\alpha,\beta) \to \R$
satisfying
\begin{equation*}
f((1-\tau)x+\tau z) \le \frac{1}{1-\tau}f(x) +\frac{1}{\tau} f(z)
\end{equation*}
whenever  $\alpha < x \le y \le z < \beta$ and $0 < \tau < 1$. Obviously, for nonnegative functions $f$, this condition is weaker
than the convexity condition \[f((1-\tau)x+\tau z) \le (1-\tau)f(x) +\tau f(z).\]
It is readily seen that, besides the nonnegative convex functions, the class $Q$ also contains all
all nonnegative monotone functions
and all nonnegative functions $f$ for which $\sup_t f(t) < 4 \inf_t f(t)$.

Two more occurrences of reciprocal Schur inequalities we want to mention are~\cite{RovTem,Watson}.
Roven\c{t}a and Temereanc\u{a}~\cite{RovTem} proved that
$f[x_1, \ldots,x_n] \ge 0$ if $f$ is a polynomial such that $f^{(n-1)} >0$, referring to~\cite{Pop}
in this connection. Paper~\cite{Watson} by G.~N.~Watson
contains a reciprocal Schur inequality in exactly our understanding.
There it is proved that if $a >0$ and $x_1, \ldots, x_n$ are distinct real numbers, then
$\sum_{j=1}^n \frac{a^{x_j}}{\prod_{k \neq j}(x_j-x_k)}$ is positive for $a >1$, of the same sign as $(-1)^{n+1}$
if $0 < a <1$, and zero for $a=1$. Watson's proof is via representing the sum as determinant and then
proceeding by induction on $n$. Obviously, the result follows from Theorem~\ref{Theo2} with
$f(t)=a^t$, in which case $f^{(n-1)}(t)=a^t(\log a)^{n-1}$.


\vspace*{1cm}

\end{document}